\documentclass[11pt]{amsart}
\usepackage{amsmath,amssymb}
\newtheorem{theorem}{Theorem}

\newtheorem{lemma}[theorem]{Lemma}

\begin{document}

\title[Minimal surfaces that pass through a prescribed point in a ball]{Area bounds for minimal surfaces that pass through a prescribed point in a ball}
\author{Simon Brendle and Pei-Ken Hung}
\address{Department of Mathematics, Columbia University, 2990 Broadway, New York, NY 10027}
\email{simon.brendle@columbia.edu}
\address{Department of Mathematics, Columbia University, 2990 Broadway, New York, NY 10027}
\email{pkhung@math.columbia.edu}
\thanks{The authors thank T\"ubingen University and the Mathematisches Forschungsinstitut Oberwolfach where part of this work was carried out. The first author was supported in part by the National Science Foundation under grant DMS-1505724.}
\begin{abstract}
Let $\Sigma$ be a $k$-dimensional minimal submanifold in the $n$-dimensional unit ball $B^n$ which passes through a point $y \in B^n$ and satisfies $\partial \Sigma \subset \partial B^n$. We show that the $k$-dimensional area of $\Sigma$ is bounded from below by $|B^k| \, (1-|y|^2)^{\frac{k}{2}}$. This settles a question left open by the work of Alexander and Osserman in 1973.
\end{abstract}
\maketitle

\section{Introduction}

In this note, we study the area of a minimal submanifold in the unit ball in $\mathbb{R}^n$. For  a minimal submanifold $\Sigma$ that passes through the center of the ball, it is well-known that the area of $\Sigma$ is bounded from below by the area of a flat $k$-dimensional disk: 

\begin{theorem}
\label{consequence.of.monotonicity}
Let $\Sigma$ be a $k$-dimensional minimal submanifold in the unit ball $B^n$ which passes through the origin and satisfies $\partial \Sigma \subset \partial B^n$. Then $|\Sigma| \geq |B^k|$.
\end{theorem}

Theorem \ref{consequence.of.monotonicity} is a direct consequence of the well-known monotonicity formula for minimal submanifolds. This technique is discussed, for example, in \cite{Ekholm-White-Wienholtz} and \cite{Simon}. 

In 1973, Alexander and Osserman \cite{Alexander-Osserman} studied a closely related problem. More precisely, they considered a minimal surface in the unit ball in $\mathbb{R}^3$ which passes through a prescribed point in the interior of the ball (not necessarily the center of the ball). In the special case of disk-type minimal surfaces, they were able to show that the area of the surface is bounded from below by the area of a flat disk. However, their argument does not work for minimal surfaces of other topological types, nor does it generalize to higher dimensions. In 1974, Alexander, Hoffman, and Osserman \cite{Alexander-Hoffman-Osserman} proved an analogous inequality in higher dimensions, but only in the special case of area-minimizing surfaces. 

In this note, we completely settle this question for minimal submanifolds of arbitrary dimension and codimension:

\begin{theorem}
\label{main.theorem}
Let $\Sigma$ be a $k$-dimensional minimal submanifold in the unit ball $B^n$ which passes through a point $y \in B^n$ and satisfies $\partial \Sigma \subset \partial B^n$. Then $|\Sigma| \geq |B^k| \, (1-|y|^2)^{\frac{k}{2}}$. Moreover, the inequality is strict unless $\Sigma$ is a flat $k$-dimensional disk which is orthogonal to $y$.
\end{theorem}

The proof of Theorem \ref{main.theorem} relies on an application of the first variation formula for minimal submanifolds (cf. \cite{Allard}, \cite{Simon}) to a carefully chosen vector field in ambient space. In particular, our argument generalizes immediately to the varifold setting. A similar technique was used in \cite{Brendle} to prove a sharp bound for the area of a free-boundary minimal surface in a ball. The main difficulty in this approach is to find the correct vector field. The vector field used in \cite{Brendle} was obtained as the gradient of the Green's function for the Neumann problem on the unit ball. By contrast, the vector field used in the proof of Theorem \ref{main.theorem} is not a gradient field, and does not have any obvious geometric interpretation.

\section{Proof of Theorem \ref{main.theorem}}

Let us fix a point $y \in B^n$. We define a vector field $W$ on $B^n \setminus \{y\}$ in the following way: For $k>2$, we define  
\begin{align*} 
W(x) 
&= -\frac{1}{k} \, \bigg ( \Big ( \frac{1-2 \langle x,y \rangle+|y|^2}{|x-y|^2} \Big )^{\frac{k}{2}} - 1 \bigg ) \, (x-y) \\ 
&+ \frac{1}{k-2} \, \bigg ( \Big ( \frac{1-2 \langle x,y \rangle+|y|^2}{|x-y|^2} \Big )^{\frac{k-2}{2}} - 1 \bigg ) \, y. 
\end{align*}
For $k=2$, we define  
\begin{align*} 
W(x) 
&= -\frac{1}{2} \, \Big ( \frac{1-2 \langle x,y \rangle+|y|^2}{|x-y|^2} - 1 \Big ) \, (x-y) \\ 
&+ \frac{1}{2} \, \log \Big ( \frac{1-2 \langle x,y \rangle+|y|^2}{|x-y|^2} \Big ) \, y. 
\end{align*}
Note that $1-2 \langle x,y \rangle+|y|^2 \geq |x-y|^2 > 0$ for all points $x \in B^n \setminus \{y\}$. This shows that $W$ is indeed a smooth vector field on $B^n \setminus \{y\}$.

\begin{lemma}
\label{a}
For every point $x \in B^n$ and every orthonormal $k$-frame $\{e_1,\hdots,e_k\} \subset \mathbb{R}^n$, we have 
\[\sum_{i=1}^k \langle D_{e_i} W,e_i \rangle \leq 1.\] 
\end{lemma}

\textbf{Proof.} 
We compute 
\begin{align*}
\sum_{i=1}^k \langle D_{e_i} W,e_i \rangle 
&= 1 - \Big ( \frac{1-2 \langle x,y \rangle+|y|^2}{|x-y|^2} \Big )^{\frac{k}{2}} \\ 
&+ \frac{(1-2 \langle x,y \rangle+|y|^2)^{\frac{k-2}{2}}}{|x-y|^k} \, \sum_{i=1}^k \langle y,e_i \rangle \, \langle x-y,e_i \rangle \\ 
&+ \frac{(1-2 \langle x,y \rangle+|y|^2)^{\frac{k}{2}}}{|x-y|^{k+2}} \, \sum_{i=1}^k \langle x-y,e_i \rangle^2 \\ 
&- \frac{(1-2 \langle x,y \rangle+|y|^2)^{\frac{k-4}{2}}}{|x-y|^{k-2}} \, \sum_{i=1}^k \langle y,e_i \rangle^2 \\ 
&- \frac{(1-2 \langle x,y \rangle+|y|^2)^{\frac{k-2}{2}}}{|x-y|^k} \, \sum_{i=1}^k \langle x-y,e_i \rangle \, \langle y,e_i \rangle \\ 
&= 1 - \frac{(1-2 \langle x,y \rangle+|y|^2)^{\frac{k}{2}}}{|x-y|^{k+2}} \, \Big ( |x-y|^2 - \sum_{i=1}^k \langle x-y,e_i \rangle^2 \Big ) \\ 
&- \frac{(1-2 \langle x,y \rangle+|y|^2)^{\frac{k-4}{2}}}{|x-y|^{k-2}} \, \sum_{i=1}^k \langle y,e_i \rangle^2 \\ 
&\leq 1.
\end{align*} 
Note that the preceding calculation is valid both for $k > 2$ and for $k=2$. This proves the assertion. \\

\begin{lemma}
\label{b}
The vector field $W$ vanishes along the boundary $\partial B^n$.
\end{lemma}

\textbf{Proof.} 
Suppose that $x \in \partial B^n$. Then $1-2 \langle x,y \rangle+|y|^2 = |x-y|^2$. This directly implies $W(x)=0$. Again, this conclusion holds both for $k>2$ and for $k=2$. \\ 

\begin{lemma}
\label{c}
We have 
\[W(x) = -(1-|y|^2)^{\frac{k}{2}} \, \frac{x-y}{k \, |x-y|^k} + o \Big ( \frac{1}{|x-y|^{k-1}} \Big )\] 
as $x \to y$.
\end{lemma} 

\textbf{Proof.} 
By definition of $W(x)$, we have 
\begin{align*} 
W(x) 
&= -(1-2 \langle x,y \rangle+|y|^2)^{\frac{k}{2}} \, \frac{x-y}{k \, |x-y|^k} + o \Big ( \frac{1}{|x-y|^{k-1}} \Big ) \\ 
&= -(1-|y|^2)^{\frac{k}{2}} \, \frac{x-y}{k \, |x-y|^k} + o \Big ( \frac{1}{|x-y|^{k-1}} \Big ) 
\end{align*} 
as $x \to y$. This proves the assertion. \\

We now describe the proof of Theorem \ref{main.theorem}. To that end, we assume that $\Sigma$ is a minimal surface in $B^n$ passing through the point $y$. Since the vector field $W$ vanishes along the boundary $\partial \Sigma \subset \partial B^n$, we obtain 
\begin{equation}
\label{first.variation}
\int_{\Sigma \setminus B_r(y)} (1-\text{\rm div}_\Sigma W) = |\Sigma \setminus B_r(y)| - \int_{\Sigma \cap \partial B_r(y)} \langle W,\nu \rangle 
\end{equation} 
by the divergence theorem. Here, $\nu$ denotes the inward pointing unit normal to the region $\Sigma \cap B_r(y)$ within the surface $\Sigma$. In other words, the vector $\nu$ is tangential to $\Sigma$, but normal to $\Sigma \cap \partial B_r(y)$. It is easy to see that 
\[\nu = -\frac{x-y}{|x-y|} + o(1)\] 
for $x \in \Sigma \cap \partial B_r(y)$. Using Lemma \ref{c}, we obtain 
\[\langle W,\nu \rangle = (1-|y|^2)^{\frac{k}{2}} \, \frac{1}{k \, r^{k-1}} + o \Big  ( \frac{1}{r^{k-1}} \Big )\] 
for $x \in \Sigma \cap \partial B_r(y)$. Since 
\[|\Sigma \cap \partial B_r(y)| = |\partial B^k| \, r^{k-1} + o(r^{k-1}),\] 
we conclude that 
\begin{equation} 
\label{term.1}
\lim_{r \to 0} \int_{\Sigma \cap \partial B_r(y)} \langle W,\nu \rangle = \frac{1}{k} \, |\partial B^k| \, (1-|y|^2)^{\frac{k}{2}} = |B^k| \, (1-|y|^2)^{\frac{k}{2}}.
\end{equation}
Combining (\ref{first.variation}) and (\ref{term.1}) gives 
\[\lim_{r \to 0} \int_{\Sigma \setminus B_r(y)} (1 - \text{\rm div}_\Sigma W) = |\Sigma| - |B^k| \, (1-|y|^2)^{\frac{k}{2}}.\] 
On the other hand, by Lemma \ref{a} we have the pointwise inequality 
\[1 - \text{\rm div}_\Sigma W \geq 0.\] 
Putting these facts together, we obtain $|\Sigma| - |B^k| \, (1-|y|^2)^{\frac{k}{2}} \geq 0$, as claimed.

Finally, we study the case of equality. Suppose that $|\Sigma| - |B^k| \, (1-|y|^2)^{\frac{k}{2}} = 0$. In this case, we have 
\[1 - \text{\rm div}_\Sigma W = 0\] 
for each point $x \in \Sigma \setminus \{y\}$. Hence, if $x$ is an arbitrary point on $\Sigma \setminus \{y\}$ and $\{e_1,\hdots,e_k\}$ is an orthonormal basis of $T_x \Sigma$, then we have 
\[|x-y|^2 - \sum_{i=1}^k \langle x-y,e_i \rangle^2 = \sum_{i=1}^k \langle y,e_i \rangle^2 = 0.\] 
This implies that $\Sigma$ is a flat $k$-dimensional disk which is orthogonal to $y$. This completes the proof of Theorem \ref{main.theorem}.

\end{document}